\documentclass[letterpaper, 10 pt, conference]{ieeeconf}  %
\IEEEoverridecommandlockouts                              %

\usepackage{graphics} %
\usepackage{epsfig} %
\usepackage{amsmath} %
\usepackage{amssymb}  %
\usepackage{booktabs}
\usepackage{makecell}
\usepackage{hyperref}

\usepackage{amsfonts}
\usepackage{xcolor}
\usepackage{booktabs}

\usepackage{optidef}
\usepackage{relsize}

\let\labelindent\relax
\usepackage{enumitem}
\usepackage{placeins}
\usepackage{comment}
\DeclareMathSymbol{\shortminus}{\mathbin}{AMSa}{"39}

\usepackage{mathtools}

\definecolor{wheat}{rgb}{0.96,0.87,0.70}

\newcommand{\ind}[1]{{_{\mathrm{#1}}}}

\newcommand{\R}{\mathbb{R}}

\newcommand{\fimpl}{f}%

\newcommand{\acados}{\texttt{acados}}

\newcommand{\code}[1]{\texttt{#1}}

\newcommand{\eye}{%
	\text{\usefont{U}{bbold}{m}{n}1}%
}

\DeclareRobustCommand{\zeros}{%
	\text{\usefont{U}{bbold}{m}{n}0}%
}

\DeclarePairedDelimiter\abs{\lvert}{\rvert}%
\DeclarePairedDelimiter\norm{\lVert}{\rVert}%

\makeatletter
\let\oldabs\abs
\def\abs{\@ifstar{\oldabs}{\oldabs*}}
\let\oldnorm\norm
\def\norm{\@ifstar{\oldnorm}{\oldnorm*}}
\makeatother

\DeclareMathSymbol{\sm}{\mathbin}{AMSa}{"39}

\newlength\figureheight
\newlength\figurewidth

\newcommand{\nsteps}{n_{\mathrm{steps}}}
\newcommand{\nstages}{n_{\mathrm{stages}}}
\newcommand{\dtotal}[2]{\frac{\mathrm{d} #1}{\mathrm{d} #2}}
\newcommand{\dpartial}[2]{\frac{\partial #1}{\partial #2}}
\newcommand{\tildeJ}{\tilde{J}}
\newcommand{\lagrange}{L}
\newcommand{\T}{{\! \top}}
\newcommand{\dtplant}{\Delta t_\mathrm{plant}}
\newcommand{\sncd}{SN}
\newcommand{\dtintegrator}{\frac{\Delta t_n}{\nsteps}}
\newcommand{\sqpstep}{d}

\newcommand{\contcontrl}{u}
\newcommand{\onlypwconstcontrol}[1]{{}}

\title{\LARGE \bf Gauss-Newton Runge-Kutta Integration for Efficient Discretization of Optimal Control Problems with Long Horizons and Least-Squares Costs}

\author{Jonathan Frey$^{1,2}$, Katrin Baumgärtner$^1$, Moritz Diehl$^{1,2}$%
\thanks{$^{1}$Department of Microsystems Engineering (IMTEK), University Freiburg, 79110 Freiburg, Germany
{\tt\small \{name.surname\}@imtek.uni-freiburg.de}}%
\thanks{$^{2}$Department of Mathematics, University Freiburg, Germany
}%
\thanks{This research was supported by DFG via Research Unit FOR 2401 and project 424107692, by BMWK via 03EI4057A and 03EN3054B, and by the EU via ELO-X 953348. }
}

\begin{document}

\maketitle
\thispagestyle{empty}
\pagestyle{empty}

\begin{abstract} %
This work proposes an efficient treatment of continuous-time optimal control problem (OCP) with long horizons and nonlinear least-squares costs.
The Gauss-Newton Runge-Kutta (GNRK) integrator is presented which provides a high-order cost integration.
Crucially, the Hessian of the cost terms required within an SQP-type algorithm is approximated with a Gauss-Newton Hessian.
Moreover, $L_2$ penalty formulations for constraints are shown to be particularly effective for optimization with GNRK.
An efficient implementation of GNRK is provided in the open-source software framework \acados.
We demonstrate the effectiveness of the proposed approach and its implementation on an illustrative example showing a reduction of relative suboptimality by a factor greater than 10 while increasing the runtime by only 10~$\%$.
\end{abstract}

\section{Introduction}

Model predictive control (MPC) is an optimization-based control strategy which relies on the (approximate) solution of nonlinear optimization problems in real-time.
In direct optimal control, starting from a continuous-time optimal control problem (OCP), a variety of choices have to be made to derive a discrete-time formulation that adequately approximates the continuous-time problem but can be solved efficiently within an online optimization context.

Direct methods for optimal control first discretize then optimize the original problem and are the focus of this paper.
Two opposing strategies %
are the single \cite{Hicks1971} and multiple shooting \cite{Bock1984} discretization approach.
While single shooting generally results in a dense NLP with less variables, multiple shooting formulations have a specific sparsity pattern and typically result in faster convergence and smaller computation times \cite{Bock1984, Albersmeyer2010}.

Sequential quadratic programming (SQP) is a widely used algorithm in the field of real-time nonlinear model predictive control (NMPC) to tackle the discretized OCP.
Especially its application via the real-time iteration (RTI) scheme \cite{Diehl2005} is of particular interest in the context of online optimization.
Within the RTI framework, a single SQP iteration is performed at each sampling time, which allows one to further split the required computation into a preparation and a feedback phase minimizing feedback delays.
Nonlinear-least squares objectives are common in control applications and enable one to use intrinsically positive-semidefinite Gauss-Newton Hessian approximations.
One essential component of SQP software for NMPC based on direct multiple shooting are integration routines that solve initial value problems with possibly nonlinear and stiff differential equations and compute the sensitivities of the result with respect to the initial state and the control input~\cite{Quirynen2017},~\cite{Frey2023}.
Often, these integration methods are simply referred to as \textit{integrators}.

The above ingredients are implemented in the open-source software package \acados{} which provides high-performance algorithms for optimal control \cite{Verschueren2021}.
It internally uses the linear algebra package \code{BLASFEO}, which provides performance-optimized routines for small to medium sized matrix operations~\cite{Frison2018}.
The \acados{} software offers a very flexible optimization problem formulation, to support a wide range of optimal-control structured problems, such as classic optimal control problems (OCP) and moving horizon estimation (MHE) problems.
Various discretization options are available, such as nonuniform grids, integrators and dedicated functionalities to handle nonlinearities and linearities in cost and constraint functions efficiently.
Moreover, a variety of quadratic programming~(QP) solution methods, such as \texttt{HPIPM}, \texttt{qpOASES}, \texttt{DAQP}, \texttt{OSQP}, \texttt{qpDUNES},
\cite{Frison2020a, Ferreau2014, Arnstrom2022, Stellato2017a, Frasch2013}
are interfaced, which either tackle the OCP-structured QP directly or after applying full or partial condensing to it~\cite{Frison2016, Axehill2015}.

This paper focuses on the discretization and Hessian approximation of the cost function.
We investigate how the control performance, both in terms of closed-loop cost and computation time, can be improved using a sophisticated cost discretization scheme.
To this end, we efficiently implemented the integration of a nonlinear least-squares Lagrange cost term together with its derivatives within an integrator, resulting in a Gauss-Newton Runge-Kutta (GNRK) integration method, recently proposed in \cite{Katliar2022}.
In addition, we propose a simple but effective penalty formulation to incorporate state constraints with an $L_2$ penalty.
The combination of the above ingredients are especially effective, in terms of accuracy and associated computational complexity, when applied to problems with relatively long horizons.
The GNRK implementation is described and its effectiveness is demonstrated together with the use of RTI and a nonuniform discretization grid in terms of computation time and closed-loop cost on the illustrative example of a pendulum on a cart.

The remainder of the paper is structured as follows.
The continuous-time OCP is presented in Section~\ref{sec:problem_formulation}.
Section~\ref{sec:discretization} discusses in detail how to transform it into an NLP using multiple shooting.
Section~\ref{sec:GNRK_algorithm} describes the GNRK integrator.
Section~\ref{sec:experiments} presents numerical experiments and Section~\ref{sec:conclusion} concludes the paper.

\section{Continuous-Time Optimal Control Problem} \label{sec:problem_formulation}
In this section, we introduce the continuous-time optimal control problem (OCP) which we aim at approximating with the direct multiple shooting formulation.

We consider optimal control problems of the form
\begin{mini!}
	{\substack{x(\cdot), \contcontrl(\cdot)}}
	{\int_{0}^{\infty} \ell(x(t), \contcontrl(t)) \, \mathrm{d} t
	\label{eq:NOCP_cost}}
	{\label{eq:NOCP}}
	{}
	\addConstraint{x(0)}{= \bar{x}_0 \label{eq:NOCP_initial_state}}
	\addConstraint{0}{= \fimpl(t, x(t), \dot{x}(t), \contcontrl(t)),~}{t\!\in\![0, \infty) \label{eq:NOCP_ODE}}
	\addConstraint{0}{\geq g(x(t), \contcontrl(t)),}{t\!\in\![0, \infty)}
\end{mini!}
where $ x(\cdot): [0, \infty) \rightarrow \R^{n_x}$, $\contcontrl(\cdot): [0, \infty) \rightarrow \R^{n_\contcontrl} $
are the state and control trajectories respectively, $ \bar{x}_0 $ is the initial state value, $ \fimpl(\cdot) $ describes the implicit system dynamics and $ g(\cdot) $ denotes the inequality constraints.
The cost function consists of the integral of
the Lagrange cost term $ \ell(\cdot) $, which we assume to have the following nonlinear least-squares form
\begin{align}
\label{eq:NLS_cost_term}
\ell(x,\contcontrl) = \frac{1}{2} \norm{r(x,\contcontrl)}^2_W,
\end{align}
where $ W \in \R^{n_y\times n_y} $ is positive definite and $r(\cdot): \R^{n_x} \times \R^{n_\contcontrl} \to \R^{n_y}$ is the potentially nonlinear residual function.

\section{Discretization of the Optimal Control Problem} \label{sec:discretization}

In this section, we want to give an overview on the various possibilities to discretize the continuous-time OCP in \eqref{eq:NOCP} within the direct multiple shooting framework.

To this end, a finite time horizon $T$ and the number of shooting intervals $N$ have to be chosen.
The shooting intervals are $[t_n, t_{n+1}]$ with $t_0 = 0, t_N = T $.
The time steps are $ \Delta t_n = t_{n+1} - t_n $ for $ n = 0, \dots, N\!\shortminus 1 $.
The multiple shooting OCP corresponding to \eqref{eq:NOCP} can then be stated as
\begin{mini!}
	{\substack{x_0,\ldots, x_N, \\ u_0,\ldots, u_{N\!\shortminus 1} }}
	{\sum_{n=0}^{N\!\shortminus 1} L_n(x_n, u_n)  + M(x_N)}
	{\label{eq:acados_OCP}}
	{}
	\addConstraint{x_0}{= \bar{x}_0}
	\addConstraint{x_{n+1}\label{eq:acados_OCP_eq}}{=\phi_n(x_n, u_n),}{n=0,\ldots,N\!\shortminus 1}
	\addConstraint{0}{\geq g_n(x_n, u_n), \label{eq:acados_OCP_ineq}}{n=0,\ldots,N\!\shortminus 1}
	\addConstraint{0}{\geq g\ind{terminal}(x_N)}.
\end{mini!}
Its optimization variables are the discrete control inputs $ u_n $ acting on $[t_n, t_{n+1}]$, $ n = 0, \dots, N\!\shortminus 1 $ and the discrete states $ x_n $ at $t_n$, $ n = 0, \dots, N $.
The values $ x_n $ and $ x_{n+1} $ are coupled by integration methods (integrators) $\phi_n(\cdot)$ that discretize the continuous-time dynamics in \eqref{eq:NOCP_ODE} and that can be different for all stages $ n=0,\dots,N\!\shortminus 1 $.
The cost terms $L_n(\cdot)$ approximate the integral of the continuous cost over the shooting interval $[t_n, t_{n+1}]$. %
The constraints $g_n(\cdot)$ represent the continuous-time constraints on $[t_n, t_{n+1}]$.
Most direct methods in optimal control only enforce the constraints at the shooting nodes.
Lastly, the terminal constraint $g\ind{terminal}(\cdot)$ and the terminal cost term $M(\cdot)$ can be used to approximately summarize the infinite remainder of the horizon.

\subsection{Constraint handling via direct penalty}
\label{sec:constraints}
In the context of NMPC, it is not recommended to impose hard constraints on the state, since this can render the OCP infeasible \cite{Rawlings2017}.
This issue is typically mitigated by \textit{softening} all constraints which depend on the state.
This means that a scalar constraint of the form $h(z) \leq 0$ in the variables $z$ is replaced by
\begin{align}
	h(z) \leq s.
	\label{eq:simple_constraint}
\end{align}
using an additional optimization variable $s$, commonly referred to as \textit{slack}, which is constrained to be nonnegative, $s \geq 0$.
The slack is penalized in the cost function, by adding a term $\rho_s(s)$, which typically consists of an $L_1$ and/or $L_2$ penalty.
Such slack variables can be considered as a control input in the context of optimization problem~\eqref{eq:acados_OCP}.
However, many OCP specific solvers allow to handle them in a more dedicated fashion, exploiting the fact that they do not enter the dynamics~\cite{Frison2020a, Frey2020}.

Since constraints are typically only imposed on the shooting nodes in the discrete-time OCP, the continuous-time trajectories corresponding to the discrete solution may violate the constraints between the shooting nodes.
This issue can be mitigated by directly adding the cost term corresponding to the constraint violation to the continuous objective, i.e.
\begin{align}
	\rho_s(\max(0, h(z))).
\end{align}
This can lead to a more accurate incorporation of the constraint cost if the cost is integrated more accurately, see below, especially if longer intervals are used.

\begin{figure}[]
	\vspace{.2cm}
	\includegraphics[width=\columnwidth]{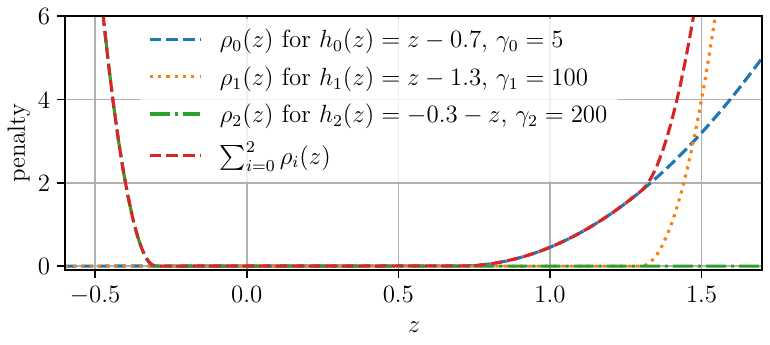}
	\caption{Multiple constraints penalized via \eqref{eq:l2_max_pen}.}
	\label{fig:max_penalties}
\end{figure}

In order to fit into the nonlinear least-squares framework, we propose to penalize constraint violations of \eqref{eq:simple_constraint} with an $L_2$ penalty with weighting parameter $\gamma$, by adding
\begin{align}
	\rho(z) = \gamma (\max(h(z), 0))^2
	\label{eq:l2_max_pen}
\end{align}
to the cost function.
The combination of such penalty functions for multiple constraints is visualized in Figure~\ref{fig:max_penalties}.
Note that for this penalty formulation, the Gauss-Newton Hessian corresponds to the exact Hessian, which is in this case not continuous.
However, the Newton iterations can be analyzed within the framework of semismooth Newton methods, compare e.g.~\cite{Hintermueller2010}.

In contrast to this, an $L_1$ penalty has a discontinuous gradient, such that a treatment as in \eqref{eq:l2_max_pen} without the square is not directly suitable for the numerical method presented in this paper.
Instead, $L_1$ penalties require the reformulation via a slack variable and inequalities.
However, an extension of formulation~\eqref{eq:l2_max_pen} to convex penalties with a continuous gradient, such as the Huber loss, would also be suitable for direct numerical treatment and can be handled with a generalized or extended Gauss-Newton Hessian (GGN / XGN)~\cite{Baumgaertner2022}.

\subsection{Cost integration scheme}
Now that state constraints are incorporated into the cost function, we want to define $L_n(\cdot)$ to approximate the continuous-time cost term on the interval $[t_n, t_{n+1})$, i.e.
\begin{align}
	L_n(x_n, u_n) \approx \int_{t_n}^{t_{n+1}}\ell(x(t), v(t)) \mathrm{d}t.
\end{align}

We consider two possible integration schemes:
\begin{enumerate}
	\item the shooting node cost discretization (\sncd), where the cost term is approximated using:
	\begin{align}
	L_n^{\mathrm{\sncd}}(x_n, u_n) & = \Delta t_n  \ell(x_n, u_n),
	\end{align}
	\item an (implicit) Runge-Kutta (RK) integration using the same integration scheme which is used to integrate the dynamics, represented by $\phi_n(\cdot)$.
	The Runge-Kutta integration of the cost and, in particular, a Hessian approximation of this cost term, which is efficient to compute, are described in detail in Section~\ref{sec:GNRK_algorithm}.
\end{enumerate}
Note that the two options coincide if the dynamics are discretized using one step of an explicit Euler scheme.

\subsection{Practical considerations and nonuniform grids}
\label{sec:togo}
In the context of real-time MPC, there are practical limitations relevant for the design of the MPC controller.
Firstly, the plant or actuators typically have a minimum sampling time $\dtplant$, for which a control input should be applied.
Operating the plant at a lower frequency restricts the control law unnecessarily, potentially sacrificing control performance.
Thus, a controller should be able to output controls with the same frequency, i.e., the sampling time of the controller $T_s$ should equal $\dtplant$, which we assume in the following.
Secondly, the control input applied for one sampling period is typically constant.
These considerations motivate a constant control input on the first shooting interval $[0, t_1]$ with $ t_1\!=\!T_s$.

While the remaining degrees of freedom in choosing a time grid are massive, the majority of practical applications of MPC restrict themselves to a uniform time grid. %
It is difficult to make general suggestions about this choice, but we want to motivate the use of a nonuniform in the following.

If the first shooting interval is longer, i.e., $\Delta t_0 > T_s $, this can lead to a loss of optimality, since the controls have to be chosen more conservatively, such that they do not drive the system away from a desired trajectory, even if applied for the longer period $\Delta t_0 $.
On the other hand, if $\Delta t_0 < T_s $, the controller might choose too aggressive actions for $T_s$, which are only safe to apply for the shorter period of $\Delta t_0$.

The prediction model and its discretization should be very accurate on the first shooting interval
to avoid suboptimality of the open-loop trajectory due to model-plant mismatch on the first part of the horizon, which is applied to the real system.

Additionally, it is essential to have a sufficiently long time horizon such that the optimizer is aware of future constraints. %
The crucial task of the latter part of the horizon which is not applied to the plant is to capture the cost-to-go accurately. %
However, since the computational cost of MPC algorithms scales at least linearly in the number of shooting nodes, a trade-off between a long time horizon and a low number of shooting nodes has to be made, which motivates the use of a nonuniform time grid.
Overall, these considerations encourage the use of a fine cost integration  in contrast to the widely used shooting node discretization.

\section{Implicit Runge-Kutta integration of the Lagrange term and Gauss-Newton Hessian Approximation}
\label{sec:GNRK_algorithm}
In the following, we describe how the Lagrange term within the continuous cost \eqref{eq:NOCP_cost} can be integrated with the same integration scheme used for the system dynamics \eqref{eq:NOCP_ODE}.
In particular, we show how the first-order, as well as an approximation of the second-order derivatives of the integrated Lagrange term can be computed with little computational overhead and a low memory footprint.
Our presentation closely follows the one given in~\cite{Katliar2022}.

\subsection{Integrated cost and its derivatives}
\label{sec:GNRK_derivation}
We assume that each integration interval $ [t_n, t_{n+1})$ is subdivided into $\nsteps$ equidistant subintervals $ [t_{n}^i, t_n^{i+1})$ with $t_n = t_{n}^0,\; t_{n+1} = t_n^{\nsteps}$ and $t_{n}^i = t_n + i \dtintegrator$ where $\dtintegrator$ is the length of each subinterval, $i=1,\dots, \nsteps$.
On each subinterval $ [t_n^i, t_n^{i+1}] $, the following system of equations is solved:
\begin{align}
s_n^{i, j} & = x_n^i + \dtintegrator \sum_{l=1}^{\nstages} a_{j,l} k_n^{i, l}, \label{eq:snij}\\
0 & = \fimpl(t_n^{i,j}, s_n^{i,j}, k_n^{i,j}, u_n),
\end{align}
for $ j=1, \ldots, \nstages$ and where $t_{n}^{i, j} = t_n^i + c_j \dtintegrator$.
The final state at the end of the subinterval is obtained as
\begin{align}
x_n^{i+1} = x_n^i  + \dtintegrator \sum_{j=1}^{\nstages} b_{j} k_{n}^{i, j}.
\end{align}

The coefficients $a_{j,l}$, $b_j$, $c_j$ are given by the Butcher tableau defining a specific RK method.
We obtain the integrated value of the cost at $t_n^i$ as:
\begin{align}
\label{eq:integrated_lagrange}
\lagrange_n^{i+1} = \lagrange_n^i + \dtintegrator  \sum_{j=1}^{\nstages} \frac{b_{j}}{2} \left\Vert r_n^{i,j} \right\Vert_W^2
\end{align}
where
$r_n^{i,j} = r(s_n^{i,j}, u_n)$.
Differentiating \eqref{eq:integrated_lagrange} with respect to $w_n := (x_n, u_n)$, we obtain
\begin{align}
	\label{integrated_lagrange_gradient}
\dtotal{\lagrange_n^{i+1}}{w_n} = \dtotal{\lagrange_n^i}{w_n} + \dtintegrator  \sum_{j=1}^{\nstages} b_{j}  r_n^{i,j^{\mathlarger{\T}}} W \tildeJ_n^{i,j},
\end{align}
where $ \tildeJ_n^{i,j} = \dtotal{ r_n^{i,j}}{w_n}$.
For the Hessian, we differentiate the above again to obtain
\begin{align} \label{eq:hessian}
\dtotal{^2\lagrange_n^{i+1}}{w_n^2} &= \dtotal{^2\lagrange_n^i}{w_n^2}  \\
&\hspace{-.3cm}\!+\!\dtintegrator \!\! \! \sum_{j=1}^{\nstages} \!\!b_{j} \!\cdot \!  \left(\tildeJ_n^{i,j^{\mathlarger{\T}}} \!W \tildeJ_n^{i,j} \!+
\!\sum_{l=1}^{n_y} r_n^{i,j,l} W_{\! l, \scriptscriptstyle{\bullet}} \dtotal{^2 r_n^{i,j,l}}{w_n^2} \right)\! \nonumber
\end{align}
where $r_n^{i,j,l}$ is the $l$-th component of $r_n^{i,j}$.
Discarding the second term within the sum in \eqref{eq:hessian}, we obtain the Gauss-Newton (GN) Hessian approximation %
\begin{align}
\dtotal{^2\lagrange_n^{i+1}}{w_n^2} \approx H_n^{i+1} := H_n^{i}
+
\dtintegrator \!\sum_{j=1}^{\nstages}\! b_{j} \tildeJ_n^{i,j^{\mathlarger{\T}}} \!W \tildeJ_n^{i,j}.
\label{eq:gnrk_hess}
\end{align}
Note that the product $ \tildeJ_n^{i,j^{\mathlarger{\T}}} \!W \tildeJ_n^{i,j} $ in \eqref{eq:gnrk_hess} is positive semidefinite.
Thus, if the coefficients $ b_j $ are nonnegative%
\footnote{Note that this is the case for Gauss-Radau IIA and Gauss-Legendre tableaus with $\nstages =1,\dots, 9$ and all explicit tableaus implemented in \texttt{acados} at time of writing, but not true in general, e.g. some DIRK methods \cite{Hairer1991} use negative $b_j$.},
the GN Hessian approximation is positive semidefinite as well.

In summary, the GNRK cost integration technique is defined by using the cost
\begin{align}
	L_n^{\mathrm{GNRK}}(x_n, u_n) & := L_n^{\nsteps}.
\end{align}
where $L_n^{\nsteps}$ as in \eqref{eq:integrated_lagrange}.
It is used together with its exact gradient \eqref{integrated_lagrange_gradient} and the GN Hessian approximation in \eqref{eq:gnrk_hess}.

If used within an SQP-type algorithm, the  Gauss-Newton Hessian approximation yields in general local linear convergence~\cite{Bock1983}.
The asymptotic linear rate depends on the deviation of the Gauss-Newton Hessian approximation from the exact Hessian at the solution~\cite{Messerer2021a}.
The deviation is explicitly given by
\begin{align} \label{eq:hessianerror}
E_n^{i+1} \!:= \! E_n^{i}
 + \dtintegrator \! \!\sum_{j=1}^{\nstages} \!b_{j}\! \left(\sum_{l=1}^{n_y} r_n^{i,j,l} W_{\! l, \scriptscriptstyle{\bullet}} \dtotal{^2 r_n^{i,j,l}}{w_n^2} \right)\!.
\end{align}
Thus, we expect fast linear convergence if the residuals $r_n^{i,j,l}$ are close to zero, i.e., the unstable upright position.

\subsection{Efficient computation of the Gauss-Newton Hessian}
\label{sec:GNRK_implementation}
This section describes how the Gauss-Newton Hessian of the integrated cost can be obtained with minor additional computations within an integrator that already delivers first-order derivatives.
Applying the chain rule, we can express $\tildeJ_n^{i,j}$ as
\begin{align}
\tildeJ_n^{i,j} = J_n^{i,j} S_n^{i,j} \label{eq:gnrk_matrix_mult}
\end{align}
where $J_n^{i,j}  = \dpartial{ r_n^{i,j}}{w} $ and
\begin{align*}
S_n^{i,j} = \dtotal{(s_n^{i,j}, u_n)}{w_n} =
\begin{bmatrix}
\dtotal{s_n^{i,j}}{x_n} & \dtotal{s_n^{i,j}}{u_n} \\
\zeros_{n_u \times n_x} & \eye_{n_u \times n_u}
\end{bmatrix}.
\end{align*}
Differentiating \eqref{eq:snij}, we obtain
\begin{align}
 \dtotal{s_n^{i,j}}{w_n} =  \dtotal{x_n^i}{w_n} + \dtintegrator \sum_{l=1}^{\nstages} a_{j, l}\dtotal{k_n^{i,l}}{w_n}.
 \label{eq:forward_sens_irk}
\end{align}

The derivatives $  \dtotal{k_n^{i,l}}{w_n}$ are available within the forward propagation of any Runge-Kutta integrator that applies internal numerical differentiation and can be reused.
For a detailed description of forward sensitivity propagation within implicit integrators, we refer to \cite{Quirynen2017}.

The additional computations for using GNRK instead of \sncd{} on a shooting interval are $\nstages\nsteps$ evaluations of $r(\cdot)$ and $\dpartial{r}{w}(\cdot)$ and matrix-matrix multiplications %
with dimension $n_y \times (n_x + n_u)$, $(n_x + n_u) \times n_x$, respectively $(n_x + n_u) \times (n_x+n_u)$.
The additional linear algebra operations are all of order $(n_x+n_u)^2$ or less, assuming that $n_y \leq (n_x+n_u)$.
Thus, it is of a lower order compared to the fastest QP solution algorithms, which require computations with order $n_x^3$ and~$n_u^3$~\cite{Frison2015a, Frey2020}.
On the other hand, in the SN discretization, the functions $r(\cdot)$ and $\dpartial{r}{w}(\cdot)$ would only be evaluated once instead of $\nstages\nsteps$ times, which might dominate the computational cost.
However, if these evaluations would dominate, the functions are likely to be very nonlinear and a more accurate integration is desirable.
Since the Hessian contributions of each point within the RK integrator can be accumulated on the fly using \eqref{eq:gnrk_hess}, \eqref{eq:gnrk_matrix_mult} and \eqref{eq:forward_sens_irk}, the additional memory footprint is small and independent of $\nsteps$ and $\nstages$.

\subsection{Comparison to generic quadrature states}
Note that a common alternative to the approach described in the previous sections is to propagate the cost with the same accuracy as the dynamics via a cost state.
Efficient integrators support a dedicated treatment of quadrature variables, i.e., variables which do not enter the DAE, such that the implicit system of equations can be decoupled and solved in the original space \cite{Hindmarsh2005}.
This idea has been extended to integrators that exploit more general linear structures within the dynamic system \cite{Quirynen2013a, Frey2019}.

The generic quadrature state approach does not take the cost function's nonlinear least-squares structure into account.
Thus, it is limited to an exact Hessian propagation, which might result in an indefinite Hessian and the associated problems within an NLP solver.
From a computational perspective, an exact Hessian propagation requires at least an additional adjoint sweep and thus more computational resources.
Different methods for Hessian propagation exist which trade-off computations and memory footprint~\cite{Quirynen2017}.

\subsection{Implementation in \acados}
The GNRK algorithm described in Section~\ref{sec:GNRK_derivation} and~\ref{sec:GNRK_implementation} has been efficiently implemented in the open-source software package \acados, which provides high-performance, embedded solvers for nonlinear optimal control.
The cost integration was implemented as an option in the \acados{} IRK module for cost functions of the nonlinear-least squares form~\eqref{eq:NLS_cost_term}.

\section{Numerical experiments}
\label{sec:experiments}
In this section, we illustrate the effectiveness of the presented strategies in terms of closed-loop cost with two numerical simulation studies.
All experiments have been carried out using \href{https://github.com/acados/acados/releases/tag/v0.2.5}{\acados{} v0.2.5} via its \texttt{Python} interface on a Laptop with an Intel i5-8365U CPU, 16 GB of RAM running Ubuntu 22.04.
The code to reproduce the results is publicly available\footnote{\url{https://github.com/FreyJo/GNRK_benchmark}}.

\begin{table}[t]
	\centering
	\vspace{.2cm}
	\caption{Overview on discretization and solver options varied in the benchmark of this paper.}
	\label{table_choices}
	\begin{tabular}{lcc}
		\toprule
		Option & Variant I & Variant II \\ \midrule
		Hessian (approximation) & Gauss-Newton (GN) & Exact Hessian (EH)                      \\
		cost discretization     & Shooting Node (SN) &Runge-Kutta (RK) \\
		discretization grid & uniform (a) & nonuniform (b) \\
		algorithm type          & converged SQP & RTI \\
		shooting intervals $N$  & 20 & 200 \\
		Time horizon $T$  & 0.4 & 4.0
		\\ \bottomrule
		\end{tabular}
\end{table}

\begin{table*}[t]
	\centering
	\vspace{.2cm}
	\caption{Closed-loop performance of different controllers measured by relative suboptimality, maximum and minimum computation time, and SQP iterations $n\ind{iter}$ over the scenario depicted in Figure~\ref{pendulum_GNRK_trajectories}.
	}
	\label{table_pendulum}
	\begin{tabular}{cccccrrrrrc}
\toprule
\begin{tabular}{@{}c@{}}Hessian approximation \\ and cost discretization \end{tabular}& $N$ & $T [\mathrm{s}]$ & RTI & uniform & rel. subopt. & max $n\ind{iter}$ & median $n\ind{iter}$ & $t_{\min} [\mathrm{ms}]$ & $t_{\max} [\mathrm{ms}]$ & in Fig.~\ref{pendulum_GNRK_trajectories}\\ \midrule
GNRK & 200 & 4.0 &  & x & \textbf{0.0} \% & 21 &4.89 &24.1& 196.1& A\\
GNSN & 200 & 4.0 &  & x & \textbf{0.0} \% & 20 &4.97 &24.5& 196.3&  \\
\midrule GNRK & 20 & 4.0 & x &  & \textbf{3.6} \% & 1 &1.00 &\textbf{0.9}& 1.2& B\\
GNRK & 20 & 4.0 &  &  & \textbf{3.7} \% & 15 &4.08 &1.8& 15.6&  \\
GNSN & 20 & 4.0 & x &  & 68.3 \% & 1 &1.00 &\textbf{0.9}& \textbf{1.1}& C\\
GNSN & 20 & 4.0 &  &  & 34.3 \% & 33 &4.42 &1.7& 30.9&  \\
EHSN & 20 & 4.0 &  &  & 34.3 \% & 400 &8.48 &1.9& 404.9&  \\
EHRK & 20 & 4.0 &  &  & \textbf{3.7} \% & 50 &5.85 &5.2& 77.0&  \\
GNRK & 20 & 4.0 & x & x & 992.3 \% & 1 &1.00 &\textbf{0.9}& \textbf{1.1}& D\\
GNRK & 20 & 4.0 &  & x & 845.4 \% & 23 &5.70 &2.9& 23.4&  \\
GNSN & 20 & 4.0 & x & x & 960.9 \% & 1 &1.00 &\textbf{0.9}& \textbf{1.0}&  \\
GNSN & 20 & 4.0 &  & x & 823.4 \% & 48 &6.38 &2.6& 45.3&  \\
GNRK & 20 & 0.4 & x & x & 3524.8 \% & 1 &1.00 &\textbf{0.9}& 1.9&  \\
GNRK & 20 & 0.4 &  & x & 3497.0 \% & 400 &178.72 &3.4& 693.9&  \\
\bottomrule 
\end{tabular}
\end{table*}

\subsection{Inverted pendulum on cart problem}
In order to demonstrate the importance of cost discretization, we regard the widely studied control problem of stabilizing an inverted pendulum mounted onto a cart.
The differential state of the model is $x = [p, \theta, s, \omega]^\top$ with cart position $p$, cart velocity $s$, angle of the pendulum $\theta$ and angular velocity $\omega$.
The control input $u$ is a force acting on the cart in the horizontal plane.
The system dynamics can be found e.g. in~\cite{Verschueren2021}.
In our OCP formulation, $u$ is constrained to be in $[-40, 40]$.
The example simulation starts with an initial state $\bar x_0 = [0, \bar \theta_0, 0, 0]^\top$ with $\bar\theta_0 = \frac{\pi}{4}$.
The goal is to drive all states to zero, i.e. the unstable upright position.
We formulate the following nonlinear least squares cost consisting of quadratic costs on states and controls and an additional penalty term penalizing a position $p$ outside of $[-1, 1] =: [p_{\min}, p_{\max}]$.
\begin{align}
	l_{\mathrm{pend}}&(x, u) = x^\top Q x + u^\top R u
	+ \rho_0(x) + \rho_1(x)
	\label{pendulum_cost}
\end{align}
with penalties $\rho_i(x)$ corresponding to the inequalities $ p_{\min} - p \leq 0$ and $ p - p_{\max} \leq 0$ with $\gamma = 5\cdot 10^4$, c.f. \eqref{eq:l2_max_pen},
and where the cost weights are $Q = \mathrm{diag}(100, 10^3, 0.01, 0.01) $, $ R = 0.2 $.
The terminal cost term is set to $M\ind{pend}(x) = x^\top P x$, where $P$ is obtained as solution of the discrete algebraic Riccati equation with cost and dynamics linearized at the steady-state.

\begin{figure}[t]
	\vspace{.2cm}
	\includegraphics[width=\columnwidth]{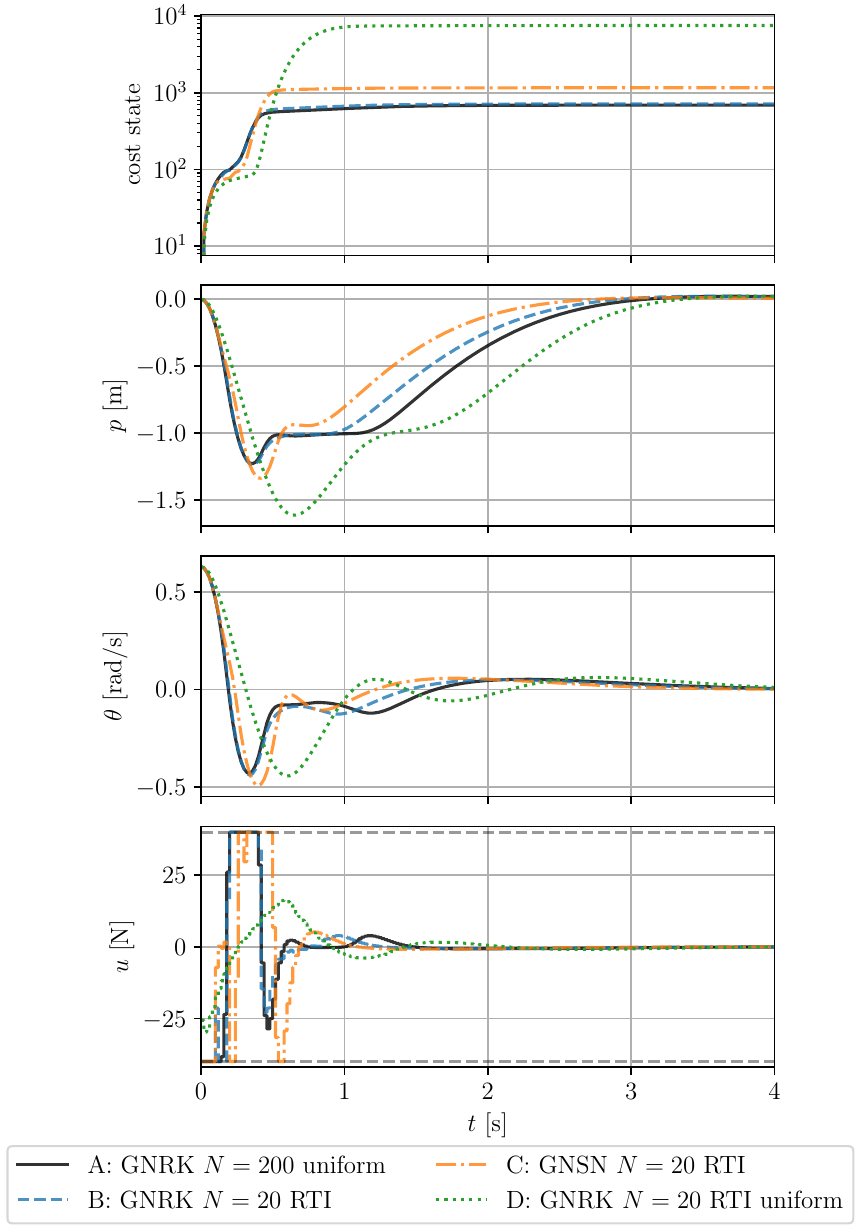}
	\caption{Closed-loop trajectories of different controllers stabilizing a pendulum on cart systems.
	The corresponding computation times are given in Table~\ref{table_pendulum}.}
	\label{pendulum_GNRK_trajectories}
\end{figure}

\subsection{Controller variants in closed-loop}
We study the behavior of different controller variants in a closed-loop simulation of $4\mathrm{s}$.
The plant is represented by a single-step IRK integrator that uses the Gauss-Radau IIA Butcher tableau with $\nstages\!=\!4 $ with a sampling time of $T_s\!=\!0.02 \mathrm{s}$.
It internally uses a model that is augmented with a cost state to accurately capture the evolution of \eqref{pendulum_cost} over time.

All controllers use \texttt{HPIPM} without condensing as a QP solver and a single step of IRK with Gauss-Radau~IIA Butcher tableau and $\nstages\!=\!4$ on each shooting interval, where the system of RK equations is solved to a tolerance of $\varepsilon_{\mathrm{IRK}}\!=\!10^{-12}$.
The time horizon of $T$ is chosen to be $4\mathrm{s}$ if not otherwise stated and is divided using one of the following time grids:
\begin{enumerate}[label=(\alph*)]
	\item uniform time grid with $\Delta t_n\!=\!\frac{T}{N}$, $n=0, \ldots, N\!\shortminus 1$
	\item nonuniform time grid using the sampling time $T_s\!=\!0.02 \mathrm{s}$ on the first interval, $\Delta t_0\!=\!T_s$, and dividing the remainder equally between the other intervals, i.e.,  $\Delta t_n\!=\!\frac{T - T_s}{N\!\shortminus 1}$, $n=1, \ldots, N\!\shortminus 1$.
\end{enumerate}
In terms of cost-discretization, we compare the shooting node (SN) and the Runge-Kutta (RK) versions, which can be combined either with the Exact Hessian (EH) or the Gauss-Newton Hessian (GN), such that the proposed approach is GNRK.
Additionally, we look at converged SQP and SQP-RTI and different number of shooting intervals $N$.
An overview on the discretization and solver options varied in this benchmark is given in Table~\ref{table_choices}.

In Figure~\ref{pendulum_GNRK_trajectories}, the closed-loop trajectories of different controllers are visualized.
Key performance indicators of even more variants are listed in Table~\ref{table_pendulum}.
The minimum and maximum computation time $t\ind{min}$, respectively $t\ind{max}$ are evaluated after running the exact same simulation 5 times and taking the minimum of each execution to remove artifacts.
The relative suboptimality is obtained by comparing the total closed-loop cost, i.e., the integrated cost state, %
and comparing it with the one of an ideal controller, i.e. without model-plant mismatch.

\begin{figure}[t]
	\centering
	\includegraphics[width=\columnwidth]{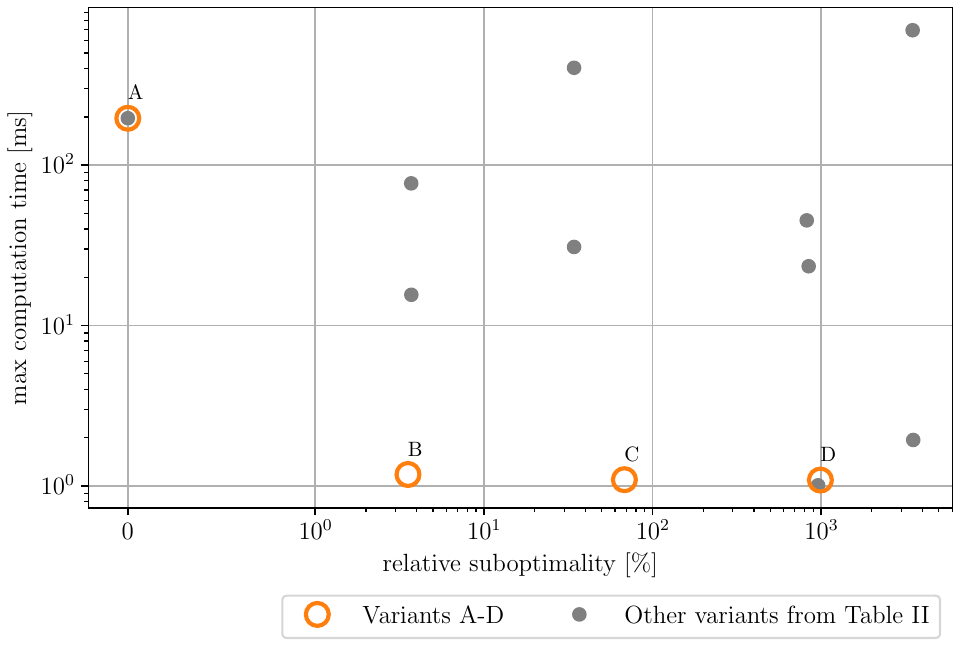}
	\caption{Pareto plot showing mean computation time and relative suboptimality of different controller variants, see Table~\ref{table_choices}. Note that the x-axis is linearly scaled in $[0, 1]$.
	}
	\label{fig:pareto_gnrk}
\end{figure}

The black line in Figure~\ref{pendulum_GNRK_trajectories} shows the reference controller with a fine uniform time grid and GNRK cost discretization.
In this case, there is no model-plant mismatch and the difference between GNRK and \sncd~cost discretization is negligible, see Table~\ref{table_pendulum}.
The controller variants with GNRK cost discretization, a nonuniform grid and $N\!=\!20$ result in a very similar closed-loop cost with only around $3.7\%$ of relative suboptimality with respect to the baseline.
In contrast, using the standard \sncd~cost discretization with otherwise the same settings, results in a suboptimality of $34\%$ and $68\%$ for converged SQP and RTI respectively.
When using controller variants with a uniform grid, with $\Delta t_0\!=\!10 T_s$, the control performance drastically degrades and we observe a relative suboptimality of over $800\%$ in this case.
On the other hand, using a uniform grid with $\Delta t_0\!=\!T_s$ results in a very short horizon length $T$, when keeping $N\!=\!20$ fix.
The corresponding controller does not staibilize the pendulum and results in over $3000\%$ of relative suboptimality.

Regarding the computation times in Table~\ref{table_pendulum}, we observe that all Gauss-Newton variants with $N\!=\!20$ and RTI have a similar runtime.
As expected, the variants with converged SQP have a much higher variance in CPU time.
Comparing the versions with exact Hessian in Table~\ref{table_pendulum} to their GN counterparts, we observe that they converge to the same solution.
However, the minimum runtime is more than twice as high.
The computation times of the GN variants with $N\!=\!200$ are roughly tenfold of the corresponding version with $N\!=\!20$.
Regarding the number of SQP iterations in Table~\ref{table_pendulum}, we observe that while the median number of iterations is similar for GNRK and GN\sncd, the maximum number is roughly half for GNRK, indicating better convergence properties.

Figure~\ref{fig:pareto_gnrk} visualizes the Pareto front of different controllers from Table~\ref{table_pendulum} in terms of relative suboptimality and maximum computation time.
The latter ultimately determines if a controller is real-time feasible.
It can be seen that the proposed controller variant B results in a reduction of relative suboptimality by a factor of 18 while increasing the maximum computation time by less than 10 \%, compared to controller variant C, i.e., the same controller without cost integration, which was regarded as an attractive variant before this work.

Overall, the results indicate that using GNRK allows one to drastically reduce the number of shooting intervals as long as the first interval is kept at $T_s$.

\subsection{Empirical contraction rate}
\begin{figure}
	\vspace{.2cm}
	\includegraphics[width=\columnwidth]{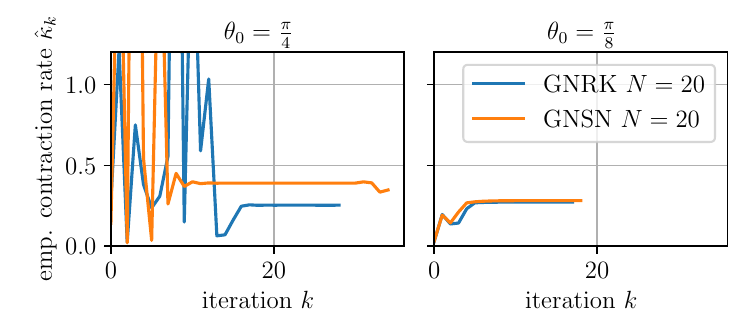}
	\caption{Empirical contraction rate for different initial states of the pendulum example. Both controllers use a uniform grid with $T=4\mathrm{s}$ and $\nstages=\!4$.}
	\label{contraction_plot}
\end{figure}
In order to verify the theoretical contraction property, we regard the empirical contraction rate
$\hat{\kappa}_k = \frac{\norm{\sqpstep_{k+1}}}{\norm{\sqpstep_{k}}},$
where $\sqpstep_k$ denotes the step in all variables at SQP iteration~$k$.
The empirical contraction rate is plotted in Figure~\ref{contraction_plot} for two initial states, $\bar x_0 = [0, \theta_0, 0, 0]^\top$.
We observe that GNRK converges with a faster rate.
For the easier initial state with $\theta_0 = \frac{\pi}{8}$, the difference is not significant, while for the initial state with $\theta_0 = \frac{\pi}{4}$, the $\hat{\kappa}_k$ values close to the solution are significantly smaller.
This gives more insight into why the maximum number of iterations in Table~\ref{table_pendulum} is higher for GNSN.

\section{Conclusion \& Outlook}
\label{sec:conclusion}
The GNRK integrator has been shown to handle nonlinear least-squares OCPs with long horizons effectively, as it trades off accuracy and computational complexity.
We showed that soft $L_2$ constraints can be handled accurately without additional slack variables by integrating the constraint violation penalty, which perfectly fits the GNRK framework.
The effectiveness of GNRK combined with the use of nonuniform discretization grids and $L_2$ penalties for state constraints has been demonstrated on an illustrative example.
This paper gave some recommendations that can help MPC practitioners to formulate and discretize their problems to obtain a competitive solver implementation.

Possible future work includes an extension of the current GNRK implementation in \texttt{acados} to handle generalized and extended Gauss-Newton Hessian terms~\cite{Baumgaertner2022}, which would allow one to handle more general convex-over-nonlinear cost and penalty functions.

\bibliography{syscop}

\end{document}